\documentclass{amsproc}
\usepackage{amssymb}

\newtheorem{theorem}{Theorem}
\newtheorem{lemma}[theorem]{Lemma}

\newtheorem{corollary}[theorem]{Corollary}

\theoremstyle{definition}

\theoremstyle{remark}
\newtheorem{remark}[theorem]{Remark}

\numberwithin{equation}{section}

\begin{document}

\title [Stability of Derivations on Hilbert $C^*$-Modules]
{Hyers--Ulam--Rassias Stability of Derivations on Hilbert
$C^*$-Modules}
\author{Maryam Amyari}
\address{Department of Mathematics, Is. Azad University, Rahnamaei Ave., Mashhad 91735, Iran, \newline and \newline Banach Mathematical Research Group (BMRG)}
\email{amyari@mshdiau.ac.ir}
\author{Mohammad Sal Moslehian}
\address{Department Mathematics, Ferdowsi University, P. O. Box 1159, Mashhad 91775, Iran, \newline and \newline Centre of Excellency in Analysis on
Algebraic Structures (CEAAS), Ferdowsi Univ., Iran.}
\email{moslehian@ferdowsi.um.ac.ir}
\thanks{The authors would like to sincerely thank
Professor Marina Haralampidou for the great hospitality during
their stay in Greece.}  \dedicatory{ Dedicated to Professor
Themistocles M. Rassias}
\begin{abstract}
Consider the functional equation ${\mathcal E}_1(f) = {\mathcal
E}_2(f)~~({\mathcal E})$~ in a certain framework. We say a
function $f_0$ is an approximate solution of $({\mathcal E})$ if
${\mathcal E}_1(f_0)$ and ${\mathcal E}_2(f_0)$ are close in some
sense. The stability problem is whether or not there is an exact
solution of $({\mathcal E})$ near $f_0$.

In this paper, the stability of derivations on Hilbert
$C^*$-modules is investigated in the spirit of
Hyers--Ulam--Rassias.

\end{abstract}

\maketitle

\section{Introduction}

One of the interesting questions in the theory of functional
equations is the following (see \cite{GRU}):
\begin{center}
{\small When is it true that a function which approximately
satisfies a functional equation ${\mathcal E}$ must be close to
an exact solution of ${\mathcal E}$?}
\end{center}
If there exists an affirmative answer we say that the equation
${\mathcal E}$ is stable.

The first stability problem was raised by S. M. Ulam during his
talk before a Mathematical Colloquium at the University of
Wisconsin in 1940 \cite{ULA}:

{\small Given a group ${\mathcal G}_1$, a metric group
$({\mathcal G}_2, d)$ and a positive number $\epsilon$, does
there exist a number $\delta>0$ such that if a function $f :
{\mathcal G}_1 \to {\mathcal G}_2$ satisfies the inequality
$d(f(xy), f(x)f(y))<\delta$ for all $x, y \in {\mathcal G}_1$
then there exists a homomorphism $T : {\mathcal G}_1 \to
{\mathcal G}_2$ such that $d(f(x), T(x))<\epsilon$ for all $x\in
{\mathcal G}_1$?}

Ulam's problem was partially solved by D. H. Hyers in 1941 in the
context of Banach spaces with $\delta = \epsilon$ in the following
form \cite{HYE}:

{\small Suppose that ${\mathcal X}_1$ and ${\mathcal X}_2$ are
Banach spaces and $f:{\mathcal X}_1 \to {\mathcal X}_2$ satisfies
the following condition: If there is $\epsilon>0$ such that
$\|f(x+y)-f(x)-f(y)\|<\epsilon$ for all $x, y \in {\mathcal
X}_1$, then there is a unique additive mapping $T:{\mathcal
X}_1\to {\mathcal X}_2$ defined by $T(x) =
\displaystyle{\lim_{n\to\infty}}\frac{f(2^nx)}{2^n}$ such that
$\|f(x)-T(x)\|<\epsilon$ for all $x\in {\mathcal X}_1$.}

Th. M. Rassias \cite{RAS1} extended Hyers' theorem in the
following form where Cauchy difference is allowed to be unbounded:

{\small Assume that ${\mathcal X}_1$ and ${\mathcal X}_2$ are
real normed spaces with ${\mathcal X}_2$ complete, $f: {\mathcal
X}_1\to {\mathcal X}_2$ is a mapping such that for each fixed
$x\in {\mathcal X}_1$ the mapping $t\mapsto f(tx)$ is continuous
on ${\mathbb R}$, and let there exist $\varepsilon\ge 0$ and
$p\in[0,1)$ such that
\[
\|f(x+y)-f(x)-f(y)\|\le \varepsilon(\|x\|^p+\|y\|^p)
\]
for all $x, y \in {\mathcal X}_1$. Then there exists a unique
linear mapping $T: {\mathcal X}_1 \to {\mathcal X}_2$ such that
\[
\|f(x)-T(x)\|\le\frac{\epsilon \|x\|^p}{1-2^{p - 1}}
\]
for all $x \in {\mathcal X}_1$.}

This result is still valid in the case where $p < 0$ by the same
approach given in \cite{RAS1} if we assume that $\|0\|^p =
\infty$. In 1990, Th. M. Rassias during the 27th International
Symposium on Functional Equations asked the question whether his
theorem can be proved for $p\geq 1$. In 1991, Z. Gajda \cite{GAJ}
following the same approach as in \cite{RAS1} provided an
affirmative solution to this question for $p>1$. Using Hyers'
method, indeed, $T(x)$ is defined by
$\lim_{n\to\infty}2\sp{-n}f(2^nx)$ if $p<1$, and
$\lim_{n\to\infty}2^nf(2\sp{-n}x)$ if $p>1$. It is shown that
there is no analogue of Th. M. Rassias' result for $p = 1$ (see
\cite{GAJ} and \cite{R-S}). This phenomenon of stability that was
introduced by Th. M. Rassias \cite{RAS1} is called {\emph
Hyers--Ulam--Rassias stability}. Thus the Hyers--Ulam stability
will be regarded as a special case of the Hyers--Ulam--Rassias
stability. A number of Rassias type results related to the
stability of various functional equations are presented in \cite{
RAS2, RAS3}.

In 1992, a generalization of Rassias' theorem was obtained by G\u
avruta as follows \cite{GAV}:

{\small Suppose $({\mathcal G},+)$ is an abelian group,
${\mathcal X}$ is a Banach space and the so-called admissible
control function $\varphi:{\mathcal G}\times {\mathcal G}\to [0,
\infty)$ satisfies
$$\widetilde{\varphi}(x, y):
 = \frac{1}{2}\displaystyle{\sum_{n = 0}^\infty}2^{-n}\varphi(2^n
x,2^n y)<\infty$$ for all $x,y\in {\mathcal G}$. If $f :
{\mathcal G} \to {\mathcal X}$ is a mapping with
$$\|f(x+y)-f(x)-f(y)\|\leq \varphi(x,y)$$ for all $x, y \in {\mathcal G}$,
then there exists a unique mapping $T : {\mathcal G} \to
{\mathcal X}$ such that $T(x+y) = T(x)+T(y)$ and
$\|f(x)-T(x)\|\leq\widetilde{\varphi}(x, x)$ for all $x,y\in
{\mathcal G}$.}

There are four methods in the study of stability of functional
equations. The first method is the {\emph direct method} in which
one uses an iteration process producing the so-called {\emph
Hyers type sequences} \cite{HYE}. Another method is based on
{\emph sandwich theorems} which are generalizations of the
Hahn-Banach separation theorems; cf. \cite{PAL}. The third
technique focuses on using {\emph invariant means}; cf.
\cite{SZE}, and the foundation of the forth method is {\emph
fixed point techniques}; cf. \cite{C-R, I-R}. The reader is
referred to \cite{CZE, H-I-R, JUN} and references therein for
further information on stability.

The notion of Hilbert $C^*$-module is a generalization of the
notion of Hilbert space.

Let ${\mathcal A}$ be a $C^*$-algebra and ${\mathcal M}$ be a
linear space which is a left ${\mathcal A}$-module with a scalar
multiplication satisfying $\lambda(xa) = x(\lambda a) = (\lambda
x)a$
for $x \in {\mathcal M},a \in {\mathcal A}, \lambda \in {\mathbb C}$.\\
The space ${\mathcal M}$ is called a {\emph pre-Hilbert
${\mathcal A}$-module} or {\emph inner product ${\mathcal
A}$-module} if there exists an inner
product $\langle .,.\rangle  :{\mathcal M} \times {\mathcal M}\to {\mathcal A}$ with the following properties:\\
(i) $\langle x,x\rangle\geq 0$; and $\langle x,x\rangle = 0~~~  {\rm iff} ~~~ x = 0;\\
(ii) \langle \lambda x + y, z\rangle = \lambda \langle x, y\rangle + \langle y,z\rangle;\\
(iii) \langle ax,y\rangle = a\langle x,y\rangle;\\
(iv) \langle x,y\rangle^* = \langle y,x\rangle$.\\
${\mathcal M}$ is called a (left) {\emph Hilbert ${\mathcal
A}$-module} if it is complete with respect to the norm $\| x \| =
\| \langle x,x\rangle\| ^{1/2}$.

{\bf (i)} Every inner product space is a left Hilbert ${\mathbb
C}$-module.

{\bf (ii)} Let ${\mathcal A}$ be a $C^*$-algebra Then every
closed left ideal $I$ of ${\mathcal A}$ is a Hilbert ${\mathcal
A}$-module if one defines $\langle a,b\rangle = ab^*\;\;(a, b\in
I)$.

Assume that ${\mathcal N}$ is another Hilbert ${\mathcal
A}$-module. Recall that a mapping $T: {\mathcal M}\to {\mathcal
N}$ is said to be adjointable if there exists a mapping $S:
{\mathcal N} \to {\mathcal M}$ such that $\langle Tx, y\rangle =
\langle x, Sy\rangle$ for all $x \in {\mathcal M}, y \in
{\mathcal N}$. The mapping $S$ is denoted by $T^*$ and called the
adjoint of $T$. If $T$ is adjointable, then it is ${\mathcal
A}$-linear and automatically continuous; cf. \cite{LAN}.

Following \cite{L-X}, a linear mapping $d: {\mathcal M} \to
{\mathcal M}$ is called a {\emph derivation} on the Hilbert
$C^*$-module ${\mathcal M}$ if it satisfies the condition
$d(\langle x, y\rangle z) = \langle d(x), y\rangle z + \langle x,
d(y)\rangle z + \langle x, y\rangle d(z)$ for every $x, y, z \in
{\mathcal M}$. It is clear that every adjointable mapping $T$
satisfying $T^* = -T$ is a derivation. The converse is not true
in general. For example, let $u_0$ be a bounded linear operator
acting on a Hilbert space $H$ of dimension greater than one such
that $u_0^* = -u_0$ and $u_0$ is not an element of the center of
$B(H)$ (e.g. consider a fixed vector $\xi \in H$ and put
$u_0(\zeta) = 2{\bf i}\langle \zeta, \xi \rangle \xi, \;(\zeta
\in H)$). Obviously, the mapping $d : B(H) \to B(H)$ defined by
$d(v) = u_0v - vu_0$ is an ordinary $*$-derivation on $B(H)$ and
so it is a derivation on $B(H)$ regarded as a Hilbert
$C^*$-module over itself. If $d$ were a $B(H)$-module map, in the
sense that $d(vw) = vd(w)$ for all $v, w \in B(H)$, then $u_0(vw)
- (vw)u_0 = v(u_0w - wu_0)$ and so $u_0vw = vu_0w$ for each $v, w
\in B(H)$. In particular, $u_0v = vu_0$ for all $v \in B(H)$
which is a contradiction. Since adjointable mappings are module
map, we infer that $d$ is not adjointable.

Recently, several extended notions of derivations such as
$(\sigma-\tau)$-derivations and generalized derivations have been
treated in the Banach algebra theory (see \cite{B-V, M-M}). In
addition, the stability of these derivations in the spirit of
Hyers--Ulam--Rassias has extensively studied by many
mathematicians; see \cite{MOS1, MOS2, PAR1, PAR2}. In this paper
we establish the stability of derivations on Hilbert
$C^*$-modules. Throughout the paper, ${\mathcal M}$ denotes a
Hilbert module over a $C^*$-algebra ${\mathcal A}$.

\section{Main results}

Recently, the stability of several mappings on Hilbert
$C^*$-modules was investigated (see \cite{AMY, MOS3}). Using some
ideas from \cite{C-R} we investigate the stability of derivations
on Hilbert $C^*$-modules. Our results may be regarded as an
extension of those of \cite{PAR2} when we consider a unital
$C^*$-algebra ${\mathcal A}$ as an ${\mathcal A}$-bimodule via
its multiplication. We start our work with a known fixed point
theorem.

\begin{theorem} [The alternative of fixed point] Suppose $({\mathcal S},d)$ is a complete
generalized metric space and $J:{\mathcal S}\to {\mathcal S}$ is a
strictly contractive mapping
with the Lipschitz constant L. Then, for each given element $x\in {\mathcal S}$,either\\
(A1) $d(J^nx, J^{n+1}x) = \infty$\\
for all $n\geq 0$, or\\
(A2) There exists a natural number $n_0$ such that:\\
(A20) $d(J^nx, J^{n+1}x)<\infty$, for all $n\geq n_0$;\\
(A21) The sequence $\{J^nx\}$ is convergent to a fixed point $y^*$ of $J$;\\
(A22) $y^*$ is the unique fixed point of $J$ in the set $U = \{y\in {\mathcal S}: d(J^{n_0}x,y) < \infty\}$;\\
(A23) $d(y,y^*)\leq\frac{1}{1-L}d(y,Jy)$ for all $y\in U$.
\end{theorem}

The following lemma gives us a useful strictly contractive
mapping.

\begin{lemma}\label{lemma}
Suppose that ${\mathcal X}$ is a Banach space, $0 \leq L < 1$ and
$\lambda \geq 0$ are given numbers and $\psi(x): {\mathcal X} \to
[0, \infty)$ has the property
\begin{eqnarray*}
\psi(x)\leq \lambda L\psi(\frac{x}{\lambda}),
\end{eqnarray*}
for all $x \in {\mathcal X}$. Assume that ${\mathcal S}: =
\{g:{\mathcal X}\to {\mathcal X}: g(0) = 0\}$ and the generalized
metric $d$ on ${\mathcal S}$ is defined by
\[d(g,h) = \inf\{c\in (0,\infty): \|g(x)-h(x)\|\leq c\psi(x),
\forall x\in {\mathcal X}\}.\] Then the mapping $J:{\mathcal S}\to
{\mathcal S}$ given by $(Jg)(x): = \frac{1}{\lambda}g(\lambda x)$
is an strictly contractive mapping.
\end{lemma}
\begin{proof}
It is easy to see that $({\mathcal S},d)$ is complete. For
arbitrary elements $g, h\in {\mathcal S}$ we have
\[
\begin{split}
d(g,h)<c &\Rightarrow \|g(x)-h(x)\|\leq c\psi(x),~~~x\in {\mathcal X}\\
&\Rightarrow \|\frac{1}{\lambda}g(\lambda
x)-\frac{1}{\lambda}h(\lambda
x)\|\leq\frac{1}{\lambda}c\psi(\lambda x),~~~x\in
{\mathcal X}\\
&\Rightarrow \|\frac{1}{\lambda}g(\lambda x)-\frac{1}{\lambda}h(\lambda x)\|
\leq Lc\psi(x),~~~x\in {\mathcal X}\\
&\Rightarrow d(Jg,Jh)\leq Lc.
\end{split}
\]

Therefore
\[
d(Jg,Jh)\leq Ld(g,h),~~~g,h\in {\mathcal S}.
\]
Hence $J$ is a strictly contractive mapping on ${\mathcal S}$
with the Lipschitz constant $L$.
\end{proof}

\begin{theorem}\label{main1} Let $\varphi:{\mathcal
M}^5 \to [0,\infty)$ be a control function with the property
\begin{equation*}
\lim_{n\to\infty}\frac{\varphi(2^nx, 2^ny, 2^nu, 2^nv, 2^nw
)}{2^n} = 0,
\end{equation*}
for all $x, y, u, v, w \in {\mathcal M}$. Suppose that
$f:{\mathcal M}\to {\mathcal M}$ is a mapping satisfying $f(0) =
0$ and
\begin{multline}\label{f}
\|f(\mu x+y)-\mu f(x)-f(y) + f(\langle u, v\rangle w) -\\
\langle f(u), v\rangle w - \langle u, f(v)\rangle w - \langle u,
v\rangle f(w)\|\leq \varphi(x, y, u, v, w),
\end{multline}
for all $\mu \in {\mathbb T} = \{z\in {\mathbb C}: |z| = 1\}$ and
$x, y, u, v, w \in {\mathcal M}$. Assume that there exists $0
\leq L_0 <1$ such that the mapping $\psi(x) =
\varphi(\frac{x}{2}, \frac{x}{2}, 0, 0, 0)$ has the property
\begin{eqnarray}\label{property1}
\psi(x)\leq 2L_0\psi(\frac{x}{2}),
\end{eqnarray}
for all $x \in {\mathcal M}$.\\
Then there exists a unique derivation $T_0:{\mathcal M}\to
{\mathcal M}$ such that
\begin{eqnarray*}
\|f(x)-T_0(x)\|\leq \frac{L_0}{1-L_0}\psi(x),
\end{eqnarray*}
for all $x\in {\mathcal M}$.
\end{theorem}

\begin{proof}
Setting $\mu = 1, y = x$ and $u = v = w = 0$ in (\ref{f}) we
obtain
\begin{equation}\label{2f}
\|f(2x)- 2f(x)\|\leq \varphi(x, x, 0, 0, 0).
\end{equation}
It follows from (\ref{property1}) and (\ref{2f}) that
\[\|\frac{1}{2}f(2x)-f(x)\|\leq \frac{1}{2}\psi(2x) \leq L_0\psi(x).\]
Hence $d(f,J_0f)\leq L_0 < \infty$, where the mapping $J_0$ is
defined on ${\mathcal S}: = \{g:{\mathcal M}\to {\mathcal M}:
g(0) = 0\}$ by $(J_0g)(x): = \frac{1}{2}g(2x)$ as in Lemma
\ref{lemma}. Applying the fixed point alternative we deduce the
existence of a mapping $T_0: {\mathcal M}\to {\mathcal M}$ such
that $T_0$ is a fixed point of $J_0$ that is $T_0(2x) = 2T_0(x)$
for all $x \in {\mathcal M}$. Since $\displaystyle{\lim_{n \to
\infty}}d(J_0^nf,T_0) = 0$ we easily conclude that
\[\lim_{n\to\infty}\frac{f(2^nx)}{2^n} = T_0(x),\]
for all $x \in {\mathcal M}$. Note that the sequence
$\{\frac{f(2^nx)}{2^n}\}$ is the Hyers sequence when one use the
direct method in establishing stability.

The mapping $T_0$ is the unique fixed point of $J_0$ in the set
$U = \{g\in {\mathcal S}: d(f,g) < \infty\}$. Hence $T_0$ is the
unique fixed point of $J_0$ such that $\|f(x) - T_0(x)\|\leq K
\psi(x)$ for some $K > 0$ and for all $x \in {\mathcal M}$.
Again, by applying the fixed point alternative theorem we obtain
\[d(f,T_0)\leq\frac{1}{1-L_0}d(f,J_0f)\leq\frac{L_0}{1-L_0},\]
and so
\[\|f(x) - T_0(x)\| \leq \frac{L_0}{1-L_0}\varphi(\frac{x}{2}, \frac{x}{2}, 0, 0, 0),\]
for all $x\in {\mathcal M}$.

It follows from (\ref{f}) that
\begin{equation}\label{fphi}
\|f(\mu x + y) - \mu f(x) - f(y)\|\leq \varphi(x, y, 0, 0, 0).
\end{equation}
Let us replace $x$ and $y$ in (\ref{fphi}) by $2^nx$ and $2^ny$,
respectively, and divide the both sides by $2^n$. Passing the
limit as $n \to \infty$ we get
\begin{equation*}
T_0(\mu x+y) = \mu T_0(x)+T_0(y),
\end{equation*}
for all $\mu \in {\mathbb T}$ and all $x, y \in {\mathcal M}$.

Next, let $\lambda\in {\mathbb C}$ ($\lambda \neq 0$) and let $K$
be a natural number greater than $4|\lambda|$. Then
$|\frac{\lambda}{K}| < \frac{1}{4} < 1 - \frac{2}{3} = 1/3$. By
Theorem 1 of \cite{K-P}, there exist three numbers $\mu_1, \mu_2,
\mu_3\in {\mathbb T}$ such that $3\frac{\lambda}{K} = \mu_1 +
\mu_2 + \mu_3$. By the additivity of $T_0$ we get
$T_0(\frac{1}{3}x) = \frac{1}{3}T_0(x)$ for all $ x\in {\mathcal
M}$. Therefore,
\begin{eqnarray*}
T_0(\lambda x) &=& T_0(\frac{K}{3}\cdot 3 \cdot
\frac{\lambda}{K}x) = \frac{K}{3}T_0(3\cdot \frac{\lambda}{K}x)\\
&=& \frac{K}{3}T_0(\mu_1x+\mu_2x+\mu_3x) = \frac{K}{3}(T_0(\mu_1x)+T_0(\mu_2x)+T_0(\mu_3x))\\
&=& \frac{K}{3}(\mu_1+\mu_2+\mu_3)T_0(x)\\
&=& \lambda T_0(x),
\end{eqnarray*}
for all $x \in {\mathcal M}$. So that $T_0$ is ${\mathbb
C}$-linear.

Putting $x = y = 0$ and replacing $u, v, w$ by $2^nu, 2^nv, 2^nw$
in (\ref{f}) we get
\begin{multline*}
\|2^{-3n}f(2^{3n}\langle u, v\rangle w) -
\langle 2^{-n}f(2^nu), v\rangle w \\
- \langle u, 2^{-n}f(2^nv)\rangle w - \langle u, v\rangle
2^{-n}f(2^nw)\| \leq 2^{-3n}\varphi(0, 0, 2^nu, 2^nv, 2^nw).
\end{multline*}
Letting $n$ tend to infinity we obtain
\[T_0(\langle u, v\rangle w) = \langle T_0(u),
v\rangle w + \langle u, T_0(v)\rangle w + \langle u, v\rangle
T_0(w),\] for all $u, v, w \in {\mathcal M}$. Hence $T_0$ is a
derivation on ${\mathcal M}$.
\end{proof}

As a consequence of Theorem \ref{main1}, we show the Rassias
stability of derivations on Hilbert $C^*$-modules.

\begin{corollary}\label{c1}
Suppose that $p \in [0, 1)$, $\alpha, \beta, \gamma >0$ and
$f:{\mathcal M} \to {\mathcal M}$ is a mapping satisfying $f(0) =
0$ and
\[
\begin{split}
\|f(\mu x + y) - \mu f(x) - f(y) - f(\langle u, v\rangle w) +
\langle f(u), v\rangle w - \langle u, f(v)\rangle w - \langle u,
v\rangle f(w)\|\\
\leq \alpha + \beta(\|x\|^p + \|y\|^p + \|u\|^p + \|v\|^p +
\|w\|^p)\\
+ \gamma
\|x\|^{p/2}\;\|y\|^{p/2}\;\|u\|^{p/2}\;\|v\|^{p/2}\;\|w\|^{p/2},
\end{split}
\]
for all $\mu \in {\mathbb T}$ and all $x, y, u, v, w \in
{\mathcal M}$. Then there exists a unique derivation
$T_0:{\mathcal M}\to {\mathcal M}$ such that
\[\|f(x)-T_0(x)\|\leq \frac{\alpha + \beta 2^{1-p}\|x\|^p + \gamma
2^{-p}\|x\|^p}{2^{1-p} - 1},\]
for all $x \in {\mathcal M}$.
\end{corollary}
\begin{proof}
Put $\varphi(x, y, u, v, w)  = \alpha + \beta(\|x\|^p + \|y\|^p +
\|u\|^p + \|v\|^p + \|w\|^p) + \gamma
\|x\|^{p/2}\;\|y\|^{p/2}\;\|u\|^{p/2}\;\|v\|^{p/2}\;\|w\|^{p/2}$,
and let $L_0 = \frac{1}{2^{1-p}}$ in Theorem \ref{main1}. Then
$\psi(x) = \alpha + 2^{1-p}\beta\|x\|^p + 2^{-p}\gamma\|x\|^p$
and there exists a derivation $T_0$ with the required property.
\end{proof}

\begin{theorem}\label{main2} Let $\varphi:{\mathcal
M}^5 \to [0,\infty)$ be a control function with the property
\begin{equation*}
\lim_{n\to\infty} 2^n \varphi(2^{-n}x, 2^{-n}y, 2^{-n}u, 2^{-n}v,
2^{-n}w ) = 0,
\end{equation*}
for all $x, y, u, v, w \in {\mathcal M}$. Suppose that the mapping
$f:{\mathcal M}\to {\mathcal M}$ is a mapping satisfying $f(0) =
0$ and
\begin{multline}\label{f2}
\|f(\mu x+y)-\mu f(x)-f(y) + f(\langle u, v\rangle w) -\\
\langle f(u), v\rangle w - \langle u, f(v)\rangle w - \langle u,
v\rangle f(w)\|\leq \varphi(x, y, u, v, w),
\end{multline}
for all $\mu \in {\mathbb T}$ and $x, y, u, v, w \in {\mathcal
M}$. Assume that there exists $0 \leq L_1 < 1$ such that the
mapping $\psi(x) = \varphi(\frac{x}{2}, \frac{x}{2}, 0, 0, 0)$
has the property
\begin{eqnarray*}
\psi(x)\leq \frac{1}{2}L_1\psi(2x),
\end{eqnarray*}
for all $x \in {\mathcal M}$.\\
Then there exists a unique derivation $T_1:{\mathcal M}\to
{\mathcal M}$ such that
\begin{eqnarray*}
\|f(x)-T_1(x)\|\leq \frac{1}{1-L_1}\psi(x),
\end{eqnarray*}
for all $x\in {\mathcal M}$.
\end{theorem}

\begin{proof}
Putting $\mu = 1, y = x$ and $u = v = w = 0$ in (\ref{f2}) we get
\begin{equation}\label{2f2}
\|f(2x)- 2f(x)\|\leq \varphi(x, x, 0, 0, 0).
\end{equation}
Replacing $x$ by $\frac{x}{2}$ in (\ref{2f2}), we obtain
\[\|f(x)-2f(\frac{x}{2})\|\leq \psi(x).\]
Hence $d(f,J_1f)\leq L_1 < \infty$, where the mapping $J_1$ is
defined on ${\mathcal S}: = \{g:{\mathcal M}\to {\mathcal M}:
g(0) = 0\}$ by $(J_1g)(x): = 2g(\frac{1}{2}x)$ as in Lemma
\ref{lemma}. Applying the fixed point alternative we deduce the
existence of a mapping $T_1: {\mathcal M}\to {\mathcal M}$ such
that $T_1$ is a fixed point of $J_1$ that is $T_1(\frac{1}{2}x) =
\frac{1}{2}T_1(x)$ for all $x \in {\mathcal M}$. Since
$\displaystyle{\lim_{n \to \infty}}d(J_1^nf,T_1) = 0$ we easily
conclude that
\[\lim_{n\to\infty} 2^n f(2^{-n}x) = T_1(x),\]
for all $x \in {\mathcal M}$.

The mapping $T_1$ is the unique fixed point of $J_1$ in the set
$U = \{g\in {\mathcal S}: d(f,g) < \infty\}$. Hence $T_1$ is the
unique fixed point of $J_1$ such that $\|f(x) - T_1(x)\|\leq K
\psi(x)$ for some $K > 0$ and for all $x \in {\mathcal M}$.
Again, by applying the fixed point alternative theorem we obtain
\[d(f,T_1)\leq \frac{1}{1-L_1}d(f,J_1f)\leq \frac{1}{1-L_1},\]
and so
\[\|f(x) - T_1(x)\| \leq \frac{1}{1-L_1}\varphi(\frac{x}{2}, \frac{x}{2}, 0, 0, 0),\]
for all $x\in {\mathcal M}$. The rest is similar to the proof of
Theorem \ref{main1}.
\end{proof}

The following corollary is similar to Corollary \ref{c1} for the
case where $p > 1$.

\begin{corollary} Suppose that $p > 1$, $\beta, \gamma >0$ and
$f:{\mathcal M} \to {\mathcal M}$ is a mapping satisfying $f(0) =
0$ and
\[
\begin{split}
\|f(\mu x+y)-\mu f(x)-f(y) + f(\langle u, v\rangle w) - \langle
f(u), v\rangle w - \langle u, f(v)\rangle w - \langle u,
v\rangle f(w)\|\\
\leq \beta(\|x\|^p + \|y\|^p + \|u\|^p + \|v\|^p + \|w\|^p) +
\gamma
\|x\|^{p/2}\;\|y\|^{p/2}\;\|u\|^{p/2}\;\|v\|^{p/2}\;\|w\|^{p/2},
\end{split}
\]
for all $\mu \in {\mathbb T}$ and all $x, y, u, v, w \in
{\mathcal M}$. Then there exists a unique derivation $T_1 :
{\mathcal M} \to {\mathcal M}$ such that
\[\|f(x)-T_1(x)\|\leq \frac{(\beta + \gamma/2)\|x\|^p}{2^{p-1}-1},\]
for all $x\in {\mathcal M}$.
\end{corollary}
\begin{proof}Put $\varphi(x, y, u, v, w)  = \alpha + \beta(\|x\|^p + \|y\|^p +
\|u\|^p + \|v\|^p + \|w\|^p) + \gamma
\|x\|^{p/2}\;\|y\|^{p/2}\;\|u\|^{p/2}\;\|v\|^{p/2}\;\|w\|^{p/2}$,
and let $L_1 = \frac{1}{2^{p-1}}$ in Theorem \ref{main2}. Then
$\psi(x) = 2^{1-p}\beta\|x\|^p + 2^{-p}\gamma\|x\|^p$ and there
exists a derivation $T_1$ with the required property.
\end{proof}

\begin{remark}
The case $p = 1$ remains unsolved.
\end{remark}

\bibliographystyle{amsalpha}

\end{document}